\documentclass{gtart_h}  

\def\ifplaintex{\expandafter\ifx\csname documentclass\endcsname\relax}

\ifplaintex 
\else
\fi

\expandafter\ifx\csname beginpicture\endcsname\relax
\expandafter\ifx\csname documentclass\endcsname\relax
\input pictex \else
\input prepictex \input pictex \input postpictex \fi\fi

\def\gt{{\mathsurround=0pt\it $\cal G\mskip-2mu$eometry \&\ 
$\cal T\!\!$opology}}        

\def\gtp{{\mathsurround=0pt\it $\cal G\mskip-2mu$eometry \&\ 
$\cal T\!\!$opology $\cal P\!$ublications}}  

\def\lognumber#1{\def\thelognumber{#1}}
\def\volumenumber#1{\def\thevolumenumber{#1}}
\def\papernumber#1{\def\thepapernumber{#1}}
\def\volumeyear#1{\def\thevolumeyear{#1}}

\def\pagenumbers#1#2{\def\startpage{#1}\def\finishpage{#2}}
\def\published#1{\def\publishdate{#1}}
\def\proposed#1{\def\theproposer{#1}}
\def\seconded#1{\def\theseconders{#1}}
\def\received#1{\def\receiveddate{#1}}
\def\revised#1{\def\reviseddate{#1}}
\def\accepted#1{\def\accepteddate{#1}}
\def\asciititle#1{\def\theasciititle{#1}}

\def\asciiemail#1{\def\theasciiemail{#1}}

\long\def\asciiabstract#1{\long\def\theasciiabstract{#1}}

\let\\\par\let\thelognumber\relax
\let\thevolumenumber\relax\let\thepapernumber\relax
\let\thevolumeyear\relax\let\thesamplenumber\relax\let\startpage\relax
\let\finishpage\relax\let\publishdate\relax\let\receiveddate\relax
\let\reviseddate\relax\let\accepteddate\relax\let\theasciititle\relax
\let\theasciiauthors\relax
\let\theasciiabstract\relax
\let\theasciiemail\relax\let\theshortauthors\relax\let\theshorttitle\relax

\long\def\maketitlep{   

\count0=\startpage

\gt\hfill      
\beginpicture
\setcoordinatesystem units <0.33truein, 0.33truein> point at 2.2 0.9
\setplotsymbol ({$\cal G$})
\plotsymbolspacing=9truept
\circulararc 315 degrees from 0 1 center at 0 0
\setplotsymbol ({$\cal T$})
\circulararc 315 degrees from 1 -1 center at 1 0
\endpicture
\break
{\small\ifx\thesamplenumber\relax 
Volume \else Sample
\fi\thevolumenumber\ (\thevolumeyear)
\startpage--\finishpage\nl
Published: \publishdate}
\vglue 0.5truein plus 0.4fil minus 0.1truein

{\parskip=0pt\leftskip 0pt plus 1fil\def\\{\par\smallskip}{\ifplaintex\large
\else\Large\fi\bf\thetitle}\par\medskip}   

\vglue 0pt plus 0.1fil 

{\parskip=0pt\leftskip 0pt plus 1fil\def\\{\par}{\sc\theauthors}
\par\medskip}

\vglue 0pt plus 0.1fil 

{\small\parskip=0pt\let\newline\\
{\leftskip 0pt plus 1fil\def\\{\par}{\sl\theaddress}\par}
\expandafter\ifx\theemail\relax    
\relax\else\vglue 5pt plus 0.02fil minus 2pt\def\\{\stdspace{\rm 
and}\stdspace} 
\cl{Email:\stdspace\tt\theemail}\fi
\ifx\theurl\relax                  
\relax\else\vglue 5pt plus 0.02fil minus 2pt\def\\{\stdspace{\rm 
and}\stdspace}
\cl{URL:\stdspace\tt\theurl}\fi\par}

\vglue 7pt plus 0.3fil minus 3pt

{\bf Abstract}
\vglue 5pt plus 0.1fil minus 2pt

\theabstract

\vglue 7pt plus 0.3fil minus 3pt

{\bf AMS Classification numbers}\quad Primary:\quad \theprimaryclass

Secondary:\quad \thesecondaryclass

\vglue 5pt plus 0.3fil minus 2pt

{\bf Keywords:}\quad \thekeywords

\vglue 10pt plus 0.5fil minus 5pt

{\small  Proposed: \theproposer\hfill Received: \receiveddate\nl
Seconded: \theseconders\hfill 
\ifx\reviseddate\relax                         
Accepted: \accepteddate                        
\else
Revised: \reviseddate                          
\fi}
\eject
}       

\let\maketitlepage\maketitlep
\let\maketitle\maketitlepage

\font\phead=cmsl9 scaled 950
\font=cmsl9 scaled 1050
\font\pnum=cmbx10 scaled 913
\font\lnum=cmbx10 
\font\pfoot=cmsl9 scaled 950
\font=cmsl9 scaled 1050
\ifplaintex
\headline{\vbox to 0pt{\vskip -4.5mm\line{\small\phead\ifnum
\count0=\startpage ISSN 1364-0380 (on line)
1465-3060 (printed) \hfill {\pnum\folio}\else\ifodd\count0\def\\{ }%
\ifx\theshorttitle\relax\thetitle\else\theshorttitle\fi\hfill{\pnum\folio}
\else\def\\{ and }{\pnum\folio}\hfill\ifx\theshortauthors\relax\theauthors
\else\theshortauthors\fi\fi\fi}\vss}}
\footline{\vbox to 0pt{\vglue 0mm\line{\small\pfoot\ifnum\count0=\startpage
\copyright\ \gtp\hfill\else
\gt, Volume \thevolumenumber\ (\thevolumeyear)\hfill\fi}\vss
}}
\else
\makeatletter
\def\@oddhead{{\small\lhead\ifnum\count0=\startpage ISSN 1364-0380 (on line)
1465-3060 (printed) \hfill {\lnum\number\count0}\else\ifodd\count0
\def\\{ }\ifx\theshorttitle\relax \thetitle \else\theshorttitle\fi\hfill
{\lnum\number\count0}\else\def\\{ and }{\lnum\number\count0}
\hfill\ifx\theshortauthors\relax 
\theauthors\else\theshortauthors\fi\fi\fi}}\def\@evenhead{\@oddhead}
\def\@oddfoot{\small\lfoot\ifnum\count0=\startpage\copyright\ \gtp\hfill\else
\gt, Volume \thevolumenumber\ (\thevolumeyear)\hfill\fi}
\def\@evenfoot{\@oddfoot}
\makeatother
\fi

\newwrite\gtoutfile
\long\gdef\makeheadfile{  
{\def\\{, }\def\s{ }
\immediate\openout\gtoutfile head.xxx
\immediate\write\gtoutfile{Proxy-for: \ifx\theasciiauthors\relax
\theauthors\else\theasciiauthors\fi\s<\ifx\theasciiemail\relax\theemail\else\theasciiemail\fi>}
\immediate\write\gtoutfile{\noexpand\\}
\immediate\write\gtoutfile{Authors: \ifx\theasciiauthors\relax
\theauthors\else\theasciiauthors\fi}
{\def\\{ }\immediate\write\gtoutfile{Title: \ifx\theasciititle\relax
\thetitle\else\theasciititle\fi}}
\immediate\write\gtoutfile{Subj-class: GT or SG or MG etc}
\immediate\write\gtoutfile{MSC-class: \theprimaryclass\ifx\thesecondaryclass\relax\else, \thesecondaryclass\fi}
\immediate\write\gtoutfile{Journal-ref: Geom. Topol. \thevolumenumber
(\thevolumeyear) \startpage-\finishpage}
\immediate\write\gtoutfile{Comments: Published by Geometry and Topology at}
\immediate\write\gtoutfile{\s\s http://www.maths.warwick.ac.uk/gt/GTVol\thevolumenumber/paper\thepapernumber.abs.html}
\immediate\write\gtoutfile{\noexpand\\}
\immediate\write\gtoutfile{}
\ifx\theasciiabstract\relax
\immediate\write\gtoutfile{\theabstract}\else
\immediate\write\gtoutfile{\theasciiabstract}\fi
\immediate\write\gtoutfile{}
\immediate\write\gtoutfile{\noexpand\\}
\immediate\write\gtoutfile{}
\immediate\closeout\gtoutfile}}  

\def\maketitlepage{\maketitlep\makeheadfile}
\let\maketitle\maketitlepage

\lognumber{346}
\received{28 July 2003}
\volumenumber{9}\papernumber{7}\volumeyear{2005}
\pagenumbers{219}{246}   
\revised{3 April 2004}
\published{26 January 2005}
\accepted{09 December 2004}
\proposed{David Gabai}
\seconded{Jean-Pierre Otal, Walter Neumann}

\usepackage{amsmath,amssymb,amscd,graphicx}

\def\no{\noindent}

\numberwithin{equation}{section}
\theoremstyle{plain}
\newtheorem{theorem}[equation]{Theorem}
\newtheorem{thm}[equation]{Theorem}
\newtheorem{proposition}[equation]{Proposition}

\newtheorem{lemma}[equation]{Lemma}
\newtheorem{corollary}[equation]{Corollary}
\newtheorem{conjecture}[equation]{Conjecture}
\newtheorem{problem}[equation]{Problem}

\theoremstyle{remark}

\theoremstyle{definition}
\newtheorem{definition}[equation]{Definition}

\newcommand{\ra}{\rightarrow}
\newcommand{\restr}{|}
\newcommand{\B}{{\mathcal B}}
\newcommand{\N}{\mathbb N}

\newcommand{\C}{\mathcal C}
\newcommand{\La}{\Lambda}

\newcommand{\R}{\mathbb R}

\newcommand{\G}{{\mathcal G}}
\renewcommand{\H}{\mathbb H}
\newcommand{\h}{{\mathcal H}}

\renewcommand{\S}{\mathbb S}

\renewcommand{\P}{{\mathcal P}}
\renewcommand{\:}{\colon}

\newcommand{\acts}{\curvearrowright}

\newcommand{\diam}{\operatorname{diam}}
\newcommand{\dist}{\operatorname{dist}}

\newcommand{\hd}{\operatorname{dist_H}}

\newcommand{\im}{\operatorname{Im}}

\newcommand{\length}{\operatorname{length}}
\renewcommand{\mod}{\operatorname{Mod}}
\newcommand{\modq}{\operatorname{Mod_Q}}

\newcommand{\reldist}{\Delta}

\newcommand{\qm}{quasi-M\"obius}
\newcommand{\sph}{\S^2}
\newcommand{\sub}{\subseteq}

\newcommand{\trip}{\operatorname{Tri}}

\newcommand{\al}{\alpha}
\def\de{\delta}
\def\De{\Delta}

\def\eps{\epsilon}
\def\ga{\gamma}
\def\Ga{\Gamma}

\def\la{\lambda}
\def\La{\Lambda}

\def\for{\quad {\rm for} \quad}

\def\Om{\Omega}

\def\ra{\rightarrow}

\def\si{\sigma}
\def\Si{\Sigma}

\def\geo{\partial_{\infty}}

\def\defeq{:=}

\newcommand{\qs}{quasisymmetric}

\begin{document}

\title{Conformal dimension and Gromov hyperbolic\\groups
with $2$--sphere boundary}
\asciititle{Conformal dimension and Gromov hyperbolic groups
with 2-sphere boundary}

\authors{Mario Bonk\\Bruce Kleiner}
\address{Department of Mathematics, University of Michigan\\Ann Arbor, 
MI, 48109-1109, USA}
\asciiemail{mbonk@umich.edu, bkleiner@umich.edu}
\gtemail{\mailto{mbonk@umich.edu}{\rm\qua and\qua}\mailto{bkleiner@umich.edu}}

\begin{abstract} Suppose $G$ is a  Gromov hyperbolic group, and
$\geo G$ is \qs ally homeomorphic to an Ahlfors  $Q$--regular
metric $2$--sphere $Z$ with Ahlfors regular 
 conformal dimension $Q$.  Then $G$ acts
discretely, cocompactly, and isometrically on $\H^3$.
\end{abstract}

\asciiabstract{%
Suppose G is a Gromov hyperbolic group, and the boundary at infinity 
of G is quasisymmetrically homeomorphic to an Ahlfors Q-regular
metric 2-sphere Z with Ahlfors regular conformal dimension Q.  Then G
acts discretely, cocompactly, and isometrically on hyperbolic 3-space.}


\maketitle

\section{Introduction}
\no
According to a well-known conjecture by Cannon, for every  Gromov hyperbolic
group $G$  whose boundary at infinity $\geo G$ 
is homeomorphic to the    $2$--sphere $\sph$, there  should  exist
a discrete, cocompact,  and isometric action of $G$  
on hyperbolic $3$--space $\H^3$. 
In the present paper we establish Cannon's conjecture  under 
the additional  assumption
that the Ahlfors regular conformal dimension of $\geo G$ is realized.

\begin{theorem} 
\label{mainthm}
Let $G$ be a Gromov hyperbolic group
with boundary $\geo G$ homeomorphic to $\sph$. If the 
Ahlfors regular conformal dimension  of $\geo G$ is attained,
then  there exists  an action of   $G$ on $\H^3$ which is 
discrete,  cocompact and isometric.
\end{theorem}

\no
By  definition, the {\em Ahlfors regular conformal dimension} of a metric 
space $Z$  is the infimal Hausdorff dimension of
all Ahlfors regular
metric spaces (see Section~\ref{prelim} for the precise definition)
 \qs ally homeomorphic to $Z$. 
This notion occurs implicitly in a paper by Bourdon and Pajot
\cite[Section 0.2]{bp3} and is a variant of Pansu's {\em conformal
dimension}  for  metric  spaces 
(the conformal dimension of a metric space $Z$
is the infimal Hausdorff dimension of  {\em all}
metric spaces \qs ally homeomorphic to $Z$).

We recall that the boundary of a Gromov hyperbolic group $G$
carries a canonical family of  {\em visual metrics}; these  are
Ahlfors regular and  pairwise quasisymmetrically homeomorphic by the
identity map.  In particular, 
it is meaningful to speak about  quasisymmetric
homeomorphisms between $\geo G$ and other metric spaces.
The assumption on the Ahlfors regular conformal dimension of $\geo G$ says 
more explicitly that there is  an  Ahlfors $Q$--regular metric space  $Z$ 
quasisymmetrically homeomorphic to $\geo G$ with smallest possible 
$Q$ among all such Ahlfors  regular  spaces.
We necessarily have $Q\ge 2$, since the Hausdorff dimension of a space
cannot be smaller than its topological dimension.
The case $Q=2$ of Theorem~\ref{mainthm}
can  easily be deduced  from 
\cite[Theorem 1.1]{qparametr} or \cite[Theorem 1.1]{quasimobius}.

The converse of Theorem \ref{mainthm} is well-known: if a group
acts discretely, cocompactly and isometrically on hyperbolic
$3$--space, then its boundary is quasisymmetrically homeomorphic
to the standard $2$--sphere \cite{Pau}, which is a $2$--regular space of 
conformal dimension $2$.  So by Theorem~\ref{mainthm}, Cannon's 
conjecture is equivalent to:

\begin{conjecture}
\label{alternate}
If $G$ is a hyperbolic group with $2$--sphere boundary,
then the Ahlfors regular conformal dimension of $\geo G$ is attained.
\end{conjecture} 

\no 
We derive Theorem~\ref{mainthm} from  \cite[Theorem 1.2]{qparametr}
and a more general result about hyperbolic groups:

\begin{theorem}
\label{loewner}
Let $Z$ be an 
 Ahlfors  $Q$--regular compact metric space, $Q>1$, 
where $Q$ is the Ahlfors regular conformal dimension of $Z$.  
If $Z$ admits a uniformly \qm\ action $G\acts Z$ which
is fixed point free and for which the induced action
on the space of triples  $\trip (Z)$ is cocompact, 
then $Z$ is $Q$--Loewner.
\end{theorem}

\no
The terminology will be explained in Section~\ref{prelim}. 
The hypotheses of  this theorem  will hold,
for example, if   $Z$ is a $Q$--regular space of 
Ahlfors regular conformal dimension $Q$, where $Q>1$, and $Z$ is
quasisymmetrically
homeomorphic to the boundary of a hyperbolic group.

Another way to state the conclusion of  Theorem~\ref{loewner}  is by saying 
that $Z$ satisfies a $(1,Q)$--Poincar\'e inequality in the
sense of Heinonen and Koskela \cite{heikos}. 
They showed 
that for    a  $Q$--regular
complete metric space such a  Poincar\'e inequality  holds if and only
if the space  is $Q$--Loewner;  they also extended  many classical
results about  quasiconformal and quasisymmmetric homeomorphisms to
the setting of $Q$--regular $Q$--Loewner  spaces.

By now there is a substantial body of  literature about 
 metric spaces satisfying   Poincar\'e
inequalities; see for example 
\cite{heikos,hajlaszkoskela,semmesnovel,semmesquant,laakso,kinshan}.
These   spaces 
play a central role in Cheeger's theory of differentiability
of Lipschitz functions \cite{cheeger}, and the Bourdon--Pajot
rigidity theorem for quasi-isometries of hyperbolic buildings
\cite{bourdonpajot}.   Theorem \ref{loewner}
suggests that one might obtain more examples of these nice spaces 
by minimizing the Hausdorff dimension
of Ahlfors regular metrics on the boundary of a hyperbolic group.
 
The full strength of the 
 group action $G\acts Z$ is actually not needed in the proof of
Theorem~\ref{loewner}.  
It is sufficient to have a collection $\G$ of
uniformly \qm\ homeomorphisms which is large enough to
map any triple in $Z$ to a uniformly separated triple, and
which does not have a common fixed point.  However,
the assumption that the action $G\acts Z$ is fixed point 
free is essential.  Starting with the Ahlfors $3$--regular
metric on $\R^2$ defined by the formula
$$d((x_1,y_1),(x_2,y_2))\defeq |x_1-x_2|+|y_1-y_2|^{{1}/{2}},$$
one can construct an Ahlfors $3$--regular metric on $\sph$ admitting 
  a uniformly
\qm\ action which is transitive on the complement of a point,
and cocompact on triples.  The sphere $\sph$ equipped with  this metric has
Ahlfors regular conformal dimension $3$,  but does not satisfy 
a $(1,p)$--Poincar\'e inequality  for any $p\ge 1$.

Similar in spirit to Theorem~\ref{mainthm} is 
another immediate 
   consequence of Theorem~\ref{loewner} for convex cocompact Kleinian 
groups.

\begin{thm} 
\label{kleinian}
Suppose   $G \acts \mathbb{H}^{n+1}$ is a   convex cocompact isometric 
action of a discrete  group $G$ on hyperbolic $n$--space $\mathbb{H}^{n+1}$, 
$n\ge 1$. 
Let $\Lambda(G)\sub \mathbb{S}^n=\geo \mathbb{H}^{n+1}$
be the limit set of $G$, and assume that 
$Q> 1$, where $Q$ is the  Hausdorff dimension of $\Lambda(G)$.
If the  Ahlfors regular conformal dimension of $\Lambda(G)$ is equal to $Q$,
then
 $Q=k\in \N$ is an integer and $\Gamma$ stabilizes a totally geodesic 
subspace of $\mathbb{H}^{n+1}$ isometric to $\mathbb{H}^{k+1}$ on which 
$\Gamma$ acts cocompactly.

\end{thm} 

\no 
Note that if under the assumptions of this theorem 
 $Z=\Lambda(G)$  carries a family
of nonconstant curves with positive $Q$--modulus, then $Q$ is equal
to the Ahlfors regular conformal dimension of $Z$ \cite[Theorem\ 15.10]{hei}.
One can also  replace the condition on the 
Ahlfors regular dimension in the previous theorem by the requirement that 
$Z$ satisfies a  $(1,p)$--Poincar\'e inequality for some  $p> 1$ (see
Section~\ref{proofsofthms}
for further discussion).

We now sketch the proof of Theorem \ref{loewner}.  Let $Z$ and 
$G\acts Z$ be as in the statement of the theorem.  A key ingredient
used repeatedly in our proof is a result of Tyson \cite{tyson} that implies
that elements of $G$ preserve  $Q$--modulus to within a 
controlled factor.  Our point of 
departure is a result of Keith and Laakso \cite{keilaa}:

\begin{theorem}[Keith--Laakso]
\label{keilaa}
Let $X$ be an 
 Ahlfors $Q$--regular complete metric space, where $Q>1$ is the
Ahlfors regular conformal dimension of $X$.  Then there exists  a
weak tangent $W$ of $X$ which carries a family of nonconstant paths  with
positive $Q$--modulus.
\end{theorem}

\no This theorem can easily be derived from 
\cite[Corollary\ 1.0.2]{keilaa}. 
For the  definition of weak tangents  and related 
discussion see \cite[Section~4]{quasimobius}; see
Section~\ref{prelim} or \cite{hei} for a discussion of modulus.
In our ``self-similar'' situation we can combine Theorem~\ref{keilaa}
with results from  \cite{quasimobius}  and 
  \cite{tyson} to obtain the following corollary, 
which may be of independent interest. 

\begin{corollary}
\label{specialkeithlaakso}
Let $Z$ be an  Ahlfors  $Q$--regular compact  metric space, where
$Q> 1$ is the Ahlfors regular conformal dimension of $Z$.
If $Z$ admits a uniformly \qm\ action $G\acts Z$ 
for which the induced action on the space of triples
$\trip (Z)$ is cocompact, then there is a  family of nonconstant paths in $Z$ with
positive $Q$--modulus.
\end{corollary}

\no
As we already pointed out, every (complete) Ahlfors $Q$--regular
space carrying a family
of nontrivial paths with positive $Q$--modulus has Ahlfors regular
   conformal 
dimension $Q$; the corollary
may be viewed as a partial converse of this fact. 

The next step in the proof of Theorem~\ref{loewner}
 is to show that $Z$ satisfies a
Loewner type condition for pairs of balls: if the $Q$--modulus for a
pair of balls is small, then their relative distance is big, quantitatively.
To prove this ball-Loewner condition,
we introduce the notion of a {\em thick path}. Thick paths
correspond to points in the support of $Q$--modulus, viewed as
an outer measure on the space of (nonconstant) paths.  Using 
the dynamics of the action $G\acts Z$, we show that any two
open sets can be joined by a thick path, and this quickly leads to 
the ball-Loewner condition.  The remaining step, which is
the bulk of our argument, shows that any complete  $Q$--regular space
satisfying the ball-Loewner condition is $Q$--Loewner.
By the result of Heinonen--Koskela mentioned above, this implies
that $Z$ satisfies a $(1,Q)$--Poincar\'e inequality.

In view of  Conjecture \ref{alternate} and Theorem \ref{loewner},
it is interesting to look for spaces  whose  Ahlfors
regular conformal dimension is (or is not) attained.
There  are  now several examples  known where the Ahlfors
regular conformal dimension is actually not realized; 
see Section~\ref{nonrealized} for more discussion. 
It is particulary interesting that Bourdon and Pajot \cite{bp3}
have found 
Gromov hyperbolic groups $G$ for which  
 $\geo G$ is not
quasisymmetrically homeomorphic to an Ahlfors regular Loewner space;
so by Theorem \ref{loewner} the Ahlfors regular conformal
dimension of $\geo G$ is not attained.

Additional remarks and open problems related to the discussion in this 
introduction can 
be found 
in the final Section~\ref{problems} of the paper.

\medskip
\textbf{Acknowledgement}\qua 
M~Bonk was  supported by NSF grants  DMS-0200566 and DMS-0244421.
B~Kleiner was supported by NSF grant  DMS-0204506.

\section{Notation and preliminaries}
\label{prelim}
\no
In this section, we will fix notation and  review some basic 
definitions and facts. We  will  be rather brief, since 
by now there  is a standard reference on these subjects \cite{hei} and 
 most of the material has been discussed in greater detail in our previous
papers \cite{quasimobius,qparametr}. 

\subsection{Notation} 
If $(Z,d)$ is a metric space, 
we denote  the open and the closed ball of radius $r>0$
centered at $a\in Z$ by $B_Z(a,r)$ and $\bar B_Z(a,r)$,
respectively. We will drop the subscript $Z$ if the space $Z$ is understood.
If $B=B(a,r)$ is a ball and $\la>0$ we let $\la B:= B(a,\la r)$. 
 We use
$\diam(A)$ for the diameter 
of a set $A\sub Z$.
If $z\in Z$ and  $A,B\sub Z$, then
$\dist(z,A)$ and   $\dist(A,B)$ are  the  distances of $z$ and $A$ and
of  $A$ and $B$, respectively.
If $A\sub Z$ and $r>0$, then we  let  $N_r(A):=\{z\in Z: \dist(z, A)<r\}$.
The
{\em Hausdorff distance} between two sets $A,B\sub Z$ is defined by
$$ \hd(A,B):= \max\big\{ \sup_{a\in A} \dist(a, B),\
\sup_{b\in B} \dist(b, A)\big\}. $$
If $f\:X\to Y$ is a map between two  spaces $X$ and $Y$,
we let $\im(f):=\{f(x): x\in X\}$. If $A\sub X$, then $f\restr A$ denotes
the restriction of the map $f$ to $A$.

\subsection{Cross-ratios and \qm\ maps}

Let $(Z,d)$ be a metric space.
The {\em cross-ratio},
$[z_1,z_2,z_3,z_4]$,
of a four-tuple of  distinct points $(z_1,z_2,z_3,z_4)$ in $Z$
is the quantity
$$[z_1,z_2,z_3,z_4]:=\frac{d(z_1,z_3)d(z_2,z_4)}{d(z_1,z_4)d(z_2,z_3)}.$$

Let $\eta\:[0,\infty)\ra [0,\infty)$
be a homeomorphism, 
and let   $f\:X\ra Y$
be an injective map between metric spaces $(X,d_X)$ and $(Y,d_Y)$.
The map $f$ is an {\em $\eta$--\qm\ map} if for
every four-tuple $(x_1,x_2,x_3,x_4)$ of distinct
points in $X$, we have
$$[f(x_1),f(x_2),f(x_3),f(x_4)]\leq \eta([x_1,x_2,x_3,x_4]).$$
The map $f$ is  {\em $\eta$--\qs\ } if
$$\frac{d_Y(f(x_1),f(x_2))}{d_Y(f(x_1),f(x_3))}
\leq \eta\left (\frac{d_X(x_1,x_2)}{d_X(x_1,x_3)}\right )$$
for every triple $(x_1,x_2,x_3)$ of distinct points in $X$.

We will make repeated use  of the following lemma. We refer to
\cite[Lemma 5.1]{quasimobius} for the proof.  
\begin{lemma}
\label{fillsin}
Let  $(Z,d)$ be  a compact metric space. 
Suppose that for each    $k\in \N$ we are given 
 a ball  $B_k=B(p_k, R_k)\sub Z$, 
 distinct points $x_k^1, x_k^2, x_k^3 \in
\bar B(p_k, \la_k R_k)$
with
$$ d_Z(x_k^i, x_k^j)> \delta_k R_k
\for i,j\in \{1,2,3\},\,  i\ne j, $$
where $\la_k, \delta_k>0$,
and an  $\eta$--\qm\ homeomorphism $g_k\:Z\ra Z$ 
such that for $y^i_k:=g_k(x_k^i)$
we have
$$d_Z(y_k^i, y_k^j)>
\delta' \for i,j\in \{1,2,3\},\, i\ne j, $$
where $\eta$ and  $\delta'>0$ are  independent of $k$.
 
\begin{itemize}
\item[\rm(i)]
If  $\lim_{k\to\infty} \la_k=0$ and the sequence
$(R_k)_{k\in \N}$ is bounded, then  
$$\diam(Z\setminus g_k(B_k)) \to 0 \for k\to \infty.  $$

\item[\rm(ii)]
Suppose for  $k\in \N$ the  set
$D_k\sub B_k$ is   $(\eps_k R_k)$--dense in
$B_k$, where $\eps_k>0$.
If $\lim_{k\to\infty} \la_k=0$ and
 the sequence   $(\eps_k/\de^2_k)_{k\in \N}$ is bounded,
then
$$ \hd(g_k(D_k), Z)\to 0 \for k\to \infty. $$
\end{itemize}
\end{lemma}
 
\no
If a group $G$ acts on a compact  metric 
space $(Z,d)$ by homeomorphisms, we write $G\acts Z$ and consider 
the group elements as self-homeomorphisms of $Z$. 
We do not require that the action is effective; so  it may well happen 
that a group element different from the unit element is represented
by the identity map.
We denote by $\trip(Z)$ the space of distinct triples in $Z$.
An action $G\acts Z$ induces an action $G\acts \trip(Z)$. The action
$G\acts \trip(Z)$ is cocompact, if and only if 
every triple in $\trip(Z)$ can be mapped to a uniformly
separated triple by a group element.
More precisely, this means that  there exists $\delta>0$ such that for
every triple $(z_1,z_2,z_3)\in \trip(Z)$ there exists a group element
$g\in G$ such that the   image triple $(g(z_1), g(z_2), g(z_3))\in \trip(Z)$
is $\delta$--separated, ie, $d(g(z_i), g(z_j))\ge \delta$ for $i\ne j$.
We call the  action $G\acts Z$  {\em fixed point free} if the maps $g\in G$
have no common fixed point, ie, for each $z\in Z$ there exists $g\in G$
such that $g(z)\ne z$. The action $G\acts Z$ is called {\em uniformly
quasi-M\"obius} if there exists  a homeomorphism $\eta\:[0,\infty)\ra
[0,\infty)$ such that every $g\in G$ is an $\eta$--\qm\ homeomorphism 
of $Z$.

\begin{lemma}
\label{qmproperties}
Suppose $(Z,d)$ is a uniformly perfect compact metric space,
and $G\acts Z$ is a fixed point free uniformly \qm\ action
which is cocompact on $\trip(Z)$.  Then the action is
minimal, ie,  for all   $z,z'\in Z$ and all  $\eps>0$, there
is a group element  $g\in G$ such that $d(g(z'),z)<\eps$.
\end{lemma}
\no 
Recall that a metric space $Z$ is called  
 {\em uniformly
perfect} if there exists a constant  $\la>1$ such that $B(a,R)\setminus
B(a, R/\la)$ is nonempty whenever $a\in Z$ and $0<R\le \diam(Z)$.

\proof
Let  $z,z'\in Z$ and  $\eps>0$ be arbitrary.   Since $Z$ is uniformly perfect,
we can choose distinct points  $x^1_k,x^2_k\in Z$ for $k\in \N$
 such that the
distances $d(z,x^1_k),\,d(z,x^2_k)$, and $d(x^1_k,x^2_k)$
agree up to  a factor independent of $k$, and 
$\lim_{k\ra\infty}d(z,x^1_k)=0$.  Set $r_k\defeq d(z,x^1_k)$.
Since $G\acts \trip(Z)$ is cocompact, we 
can find  $g_k\in G$ such that the triples
$(g_k(z),g_k(x^1_k),g_k(x^2_k))$ are $\de$--separated
where $\de>0$ is independent of $k$. Choose $R_k>0$ such that 
$\lim_{k\ra\infty}R_k=0$ and $\lim_{k\ra\infty}R_k/r_k=\infty$. 
  By Lemma~\ref{fillsin} we then have 
$\lim_{k\ra\infty}\diam(Z\setminus g_k(B_k))=0$, where $B_k:= B(z,R_k)$.
Pick $g\in G$ such that $g(z')\neq z'$.  Then for large
$k$, either $z'\in g_k(B_k)$ or $g(z')\in g_k(B_k)$,
which means that either $g_k^{-1}(z')\in B_k$
or $g_k^{-1}\circ g(z')\in B_k$.  Hence  for sufficiently large
$k$ one of the points
$g_k^{-1}(z')$ or $g_k^{-1}\circ g(z')$ is contained in 
$B(z,\eps)$.
\qed

\subsection{Modulus in metric measure spaces}
Suppose  $(Z,d,\mu)$ is a metric measure  space, ie, $(Z,d)$ is a
metric space  and $\mu$ a Borel measure on $Z$. Moreover, we assume that 
that $(Z,d)$ is locally compact and that   $\mu$
is locally finite and has dense support.

The space $(Z,d,\mu)$ is called ({\em Ahlfors})
$Q$--{\em regular}, $Q > 0$, if the measure $\mu$ satisfies
\begin{equation} \label{regular}
C^{-1}  R^Q \le \mu(B(a,R)) \le CR^Q
\end{equation}
for each open ball $B(a,R)$ of radius $0 < R \le \diam(Z)$ and for some
constant $C \ge 1$ independent of the ball.  If the measure is not 
specified, then it is understood that $\mu$ is $Q$--dimensional 
Hausdorff measure. Note that a complete Ahlfors regular space $Z$ is 
uniformly perfect and  {\em 
proper}, ie, closed balls in $Z$ are compact. 

A  {\em density} (on $Z$)  is a Borel function $\rho\:Z\ra [0,\infty]$.
A density $\rho$ is called  {\em admissible} for a path family
$\Gamma$  in $Z$, if
$$      \int_\ga \rho\, ds \ge 1                $$
for each  rectifiable  path
  $\ga \in \Gamma$. Here integration is with respect to arclength
on $\ga$. 
 If $Q\ge 1$,
the $Q$--{\it modulus} of a family
$\Gamma$ of paths in $Z$ is the number
\begin{equation}
\mod_Q (\Gamma) = \inf \int \rho^Q \, d\mu ,
\end{equation}
where the infimum is taken over all densities  $\rho \: Z \ra [0,\infty]$
that are admissible for $\Gamma$.
If $E$ and $F$ are subsets of $Z$ with positive diameter,
we denote by 
\begin{equation} \label{reldist}
   \Delta(E,F) \defeq \frac {\dist(E,F)} {\min\{\diam(E), \diam(F)\}   }
\end{equation}
the {\em relative distance} of $E$ and $F$,  and by $\Ga(E,F)$ 
 the family of  paths in $Z$ connecting
$E$ and $F$.

Suppose $(Z,d,\mu)$ is a connected metric measure space.
Then $Z$ is called a $Q$--{\it Loewner space}, $Q \ge  1$,  if
there exists a positive decreasing  function $\Psi\: (0,\infty)
\ra (0,\infty)$ such that
\begin{equation} \label{defloewner}
\mod_Q(\Ga (E,F)) \ge \Psi(\Delta(E,F))
\end{equation}
whenever $E$ and $F$ are
disjoint
continua in $Z$. Note that in \cite{heikos} it was also required 
 that $Z$ is rectifiably connected. 
In case that the (locally compact) space 
 $(Z,d,\mu)$ is $Q$--regular and $Q>1$, it is unnecessary to 
make this additional assumption, because  property 
(\ref{defloewner}) alone implies that $(Z,d)$ is even  quasiconvex, ie, 
for  every pair of points there exists a connecting path whose 
length is comparable 
to the distance of the points.

We will need the following result due to  Tyson. 

\begin{theorem}
\label{qmodpreserved}
Let $X$ and $Y$ be Ahlfors $Q$--regular locally compact  metric spaces,
$Q\ge 1$, 
 and
let $f\:X\ra Y$ be an $\eta$--\qm\ homeomorphism.  Then
for every family $\Ga$ of paths in $X$, we have
$$
\frac{1}{C}\modq(\Ga)\leq\modq(f\circ\Ga)\leq C\modq(\Ga),
$$
where $f\circ\Ga\defeq\{f\circ\ga : \ga\in\Ga\}$ and
C is a constant depending only on $X$, $Y$ and $\eta$.
\end{theorem}

\no
Tyson proved this for quasisymmetric mappings $f$ in \cite{tyson} and for
locally quasisymmetric maps in \cite[Theorem 6.4 and Lemma 9.2]{tyson2}.
Here a map $f\:X\ra Y$ is called {\em locally $\eta$--\qs\ }
 if every point $x\in X$
has an open neighborhood $U$ such that $f|U$ is $\eta$--quasisymmetric.
Since $\eta$--\qm\ maps are locally
$\tilde\eta$--\qs\ with $\tilde\eta$ depending
only on $\eta$, the above theorem follows.

\begin{lemma}
\label{qmodproperties}
Assume $Q>1$ and $(Z,d,\mu)$ is Ahlfors $Q$--regular.  
Then there exists a constant $C>0$ with the following property. 
If $\Gamma$ is a family of paths in $Z$ which start in a ball $B\sub Z$
of radius $R>0$ and whose lengths are  at least $LR$, where 
$L\ge 2$, then
$$
\mod_Q(\Ga)<C\left(\log L\right)^{1-Q}. 
$$
\end{lemma}
\no
We omit a detailed proof, since the statement is well-known
(cf\ \cite[Lemma 7.18]{hei} and \cite[Lemma 3.2]{BKR} for very
similar results). 
The basic idea is to use a test function  of the form 
$$  \frac{c}{ R+\dist(x, B)}$$ supported in 
$LB$ and use the upper mass bound $\mu(B(x,r))\lesssim r^Q$.

\subsection{Thick paths}
We now assume that  $(Z,d,\mu)$  is a separable locally compact  metric
measure space, and $Q\ge 1$. Let $I:=[0,1]$, and denote by 
$\P\defeq C(I,Z)$ the set of 
(continuous) paths
in $Z$, metrized by the supremum  metric. Then   $\P$ is a 
separable complete metric space.
 Since $Q$--modulus is monotonic and countably 
subadditive for path families (cf\ \cite[p.\ 51]{hei}), we can consider 
$\mod_Q$ as an  
outer measure on $\P$.

\begin{definition}
\label{defthick}
A path  $\ga\in\P$ is {\em thick} if for all $\eps>0$, the family
of nonconstant paths  in the ball 
$B(\ga,\eps)\sub \P$ has positive $Q$--modulus.  
\end{definition}
\no
In other words, a path $\ga\in \P$ is thick if
$\mod_Q(B(\ga, \eps)\setminus\C)>0$ for all $\eps>0$,
where $\C$ is  the family of 
constants paths in $Z$.
We have to  exclude the constant paths here, because 
$\mod_Q(\Ga)=\infty$ whenever $\Ga$ contains such a path. 
Constant paths lead to some technicalities later on, which could be avoided 
if we had defined $\P$ as the space of nonconstant paths 
in $Z$. This also has disadvantages, since 
certain completeness and compactness properties 
of  $\P$ would be lost  with this definition.

We denote by  $\P_T$  the set of
thick paths in $\P$.   
Ignoring constant paths, the thick paths form the support of the
outer measure $\mod_Q$.

\begin{lemma}\textnormal{(Properties of thick paths)}
\label{thickproperties}
\begin{itemize}
\item[\rm{ (i)}] \textnormal{(Stability under limits)}\qua
The set  $\P_T$ is closed
in $\P$: if
$\ga_k\in\P$ is thick for $k\in \N$ and $\ga=\lim_{k\ra\infty}\ga_k$,
 then $\ga$ is thick.
\smallskip 

\item[\rm{ (ii)}] \textnormal{(Thickness of subpaths)}\qua
The composition of any embedding $I\ra I$
with a thick path is a thick path.  
\smallskip 

\item[\rm{ (iii)}]\textnormal{(Quasi-M\"obius invariance)}\qua
 If $(Z,d,\mu)$ is locally compact and  Ahlfors $Q$--regular, $Q\ge 1$, then
the image of a thick path under a  \qm\ homeomorphism
$Z\ra Z$ is thick. 

\end{itemize}
\end{lemma}
\proof
Property (i) follows  immediately from the definitions. 
Property (iii) is a consequence of Tyson's  Theorem~\ref{qmodpreserved}. 

To prove property (ii) first note that if $\ga$ is any path 
and $\al\:I\ra I$ is  an embedding, then the definition
of modulus implies that 
\begin{equation}
\label{restrict}
\modq(B(\ga,\eps))\leq \modq(B(\ga\circ\al,\eps)).
\end{equation}
If $\ga$ is a thick  path and $\ga\circ\al$ is nonconstant,
then for sufficiently small  $\eps>0$ there will be no constant
paths in either $B(\ga,\eps)$ or $B(\ga\circ\al,\eps)$, and
(\ref{restrict}) implies that $\ga\circ\al$ is thick.
If $\ga$ is thick and $\ga\circ\al$ is constant, we can assume
without loss of generality that $\im(\al)$ is not contained
in a larger interval on which $\ga$ is constant. If
$\im(\al)=I$ then  $\ga\circ\al$ is just
a reparametrization of a thick path and is therefore thick.
Otherwise, we can enlarge  $\im(\al)$ slightly
and approximate $\ga\circ \al$ 
by nonconstant    subpaths
of $\ga$, which are thick by the first part of the argument. 
Property (i) now  implies that $\ga\circ\al$ is also  thick.
\qed 

\begin{lemma}
\label{aethick} The family  $\P_t$ of nonconstant    paths in $\P$ which are 
not thick 
has zero $Q$--modulus.
In particular, given any   family  
$\Ga\sub \P$ of nonconstant paths,
we have $\modq(\Ga\cap\P_T)=\modq(\Ga)$. 
\end{lemma}

\proof
For each $\ga\in\P_t$, we can find an open set $U_\ga$ 
containing $\ga$ which consists of nonconstant 
paths and has zero  $Q$--modulus.  The space $\P$ is separable,
so we can find a countable subcollection of the sets   $U_\ga$
which covers $\P_t$. Countable subadditivity of 
$Q$--modulus  implies that $\modq(\P_t)=0$.

The second part of the  lemma follows from monotonicity and subadditivity 
of $Q$--modulus.
\qed 

\medskip\no
The previous lemma implies the existence of nonconstant 
 thick paths whenever $Z$  
carries a family of nonconstant paths of positive $Q$--modulus. 
Moreover, suppose 
$\Ga_0$  is a 
   family 
 of  paths with $\modq(\Ga_0)=0$. Then if $\ga$ is thick and $\eps>0$ 
is arbitrary, we can  find a thick path $\al\in B(\ga, \eps) \setminus
\Ga_0$.  In other words, by a small perturbation of a thick path, we can 
obtain a thick  path   avoiding any  
given family   of zero $Q$--modulus.

\section{The Loewner condition for balls}
\label{ballloewnersection}

\no
In this section we prove the following proposition,
which asserts that a space which satisfies a Loewner
type condition for pairs of balls, satisfies the
Loewner condition for all pairs of continua.  

\begin{proposition}
\label{ballloewner}
Let $(Z,d,\mu)$ be a complete  metric measure space.
Assume that for all $C>0$, there are constants $m=m(C)>0$
and $L=L(C)>0$ such that if $R>0$ and $B,B'\sub Z$ are
$R$--balls with $\dist(B,B')\leq CR$, then the
$Q$--modulus of the family
$$
\left\{\ga\in\Ga(B,B'):\length(\ga)\leq LR\right\}
$$
is greater than $m$.   Then $(Z,d,\mu)$ is  
$Q$--Loewner. 
\end{proposition}

\no
Rather than using the hypothesis directly, the proof of the proposition
will use the following consequence:
 if $\rho\:Z\ra [0,\infty]$ is a Borel
function and the balls are as in the statement of the proposition, 
then there is a path $\si\in\Ga(B,B')$ whose length is at most 
$LR$ and  whose $\rho$--length is at most
\begin{equation}
\label{rholengthbound}
\frac{1}{m^{{1}/{Q}}}\left(\int_{(L+1)B}\rho^Q\,d\mu\right)^{{1}/{Q}}.
\end{equation}
Here and in the following we call the integral $\int_\sigma \rho \, ds$
the {\em $\rho$--length}  of a rectifiable path $\sigma$. 

We point out the following corollary of Proposition~\ref{ballloewner} which is 
of independent interest. 

\begin{corollary} \label{ballloewnercor}
Suppose $(Z,d,\mu)$ is an Ahlfors $Q$--regular complete metric measure space,
$Q>1$. Suppose that there  exists a positive decreasing function 
$\Psi\: (0,\infty) \ra (0,\infty)$ such that 
$$ \modq(\Ga(B,B')) \ge \Psi(\Delta(B,B'))$$ 
whenever $B$ and $B'$ are disjoint balls in $Z$. 
Then $(Z,d,\mu)$ is $Q$--Loewner.  
\end{corollary}

\proof This immediately follows from Proposition~\ref{ballloewner}, 
Lemma~\ref{qmodproperties} and the subadditivity of modulus (see the end 
of the proof of Lemma~\ref{conctballs} for additional details).
\qed 

\medskip
\no
Before we start with the 
proof of Proposition~\ref{ballloewner}, 
we first indicate a lemma whose proof uses  a similar construction 
 in  a simpler setting.

\begin{lemma}
Let $(X,d)$ be  a complete metric space.  Suppose there
exist $0\leq \la<{1}/{2}$ and $L<\infty$ such that
if $u,v\in X$, then there is a path of length
at most $Ld(u,v)$ connecting  $B(u,\la d(u,v))$
and $B(v,\la d(u,v))$.  Then $X$ is quasi-convex.
\end{lemma}

\no
This lemma can be used to give another proof that a $Q$--regular
space satisfying a $(1,Q)$--Poincar\'e inequality is quasi-convex.

\medskip
\no
{\bf Outline of   proof}\qua  Suppose $x,y\in X$ and let $R:=d(x,y)$. 
By assumption we can find a path $\si_1$ of length $\leq LR$
joining $B(x,\la R)$ to $B(y,\la R)$.  Set $\Si_0\defeq\{\si_1\}$.
Then we can find paths $\si_2,\,\si_3$ of length
$\leq L\la^2R$ such that $\si_2$ joins 
$B(x,\la^2R)$ to $B(\si_1(0),\la^2R)$ and 
$\si_3$ joins $B(\si_1(1),\la^2R)$ to $B(y,\la^2R)$.
Set $\Si_2\defeq\{\si_2,\si_3\}$.  Continuing inductively, we
can find path collections $\Si_k$ for all $k\geq 0$. 
At each step of the induction  the ``total gap'' $\De_k$
gets multiplied by $2\la<1$, and the total length
of the curves generated is $\leq L\De_k$.  One then concludes
that if 
$$Y\defeq\bigcup_{k=0}^\infty\bigcup_{\si\in\Si_k}\im(\si),
$$
then $\bar Y$ is a compact connected set containing $\{x,y\}$,
and has $1$--dimensional Hausdorff measure at most
$\frac{LR}{1-2\la}$.  Therefore there is 
an arc of length 
$$\leq \left(\frac{L}{1-2\la}\right) d(x,y)
$$
joining $x$ to $y$.
\qed
\medskip

\no 
The proof of Proposition \ref{ballloewner} will require two lemmas.

\begin{lemma}
\label{goodpoints}
Let $X$ be a metric space, and $\nu$ be a 
Borel measure on $X$.  If $Y\sub X$ is a nondegenerate 
continuum, 
then we can find a point $y\in Y$ such that
for all $r>0$ we have
\begin{equation}
\label{smallmass}
\nu(B(y,r))\leq \frac{10 r }{\diam(Y)} \nu(X).
\end{equation}
\end{lemma}
\proof
We assume that $0<\nu(X)<\infty$, for otherwise the assertion obviously
holds.

If the statement were false, then for each  $y\in Y$ we could find 
$r_y>0$ such that $\nu(B(y,r_y))>M\nu(X) r_y$, where $M=10/\diam(Y)$.
Then the  radii  $r_y$, $y\in Y$, are uniformly bounded from above,
and so 
we can find a disjoint subcollection $\{B_i=B(y_i,r_i)\}_{i\in I}$
of the cover $\{B(y,r_y): y\in Y\}$ of $Y$
such that the collection $\{5B_i\}_{i\in I}$ is also a  cover of 
 $Y$ \cite[Theorem 1.2]{hei}.
Define an equivalence relation on $I$ by declaring
that $i\sim i'$ if there are elements 
$i=i_1,\ldots,i_k=i'$ such that $5B_{i_j}\cap 5B_{i_{j+1}}\neq\emptyset$
for $1\leq j< k$.  If $I=\sqcup_{j\in J}I_j$ is the 
decomposition of $I$ into equivalence classes, then the
collection $\{\cup_{i\in I_j}5B_i\}_{j\in J}$
is a cover of $Y$ by disjoint open sets; since $Y$
is connected this implies that $\#J=1$.  It follows
that 
$$
\sum_{i} r_i\geq \frac 1{10} \diam(Y) =\frac1M,
$$
and so
$$\nu(X)\geq \sum_i \nu(B_i) >  \nu(X) M  \sum_i r_i \ge \nu(X). $$
This is a contradiction. 
\qed

\begin{lemma}
\label{subballs}
Let $Z$ be as in Proposition \ref{ballloewner}, and 
suppose $0<\la<{1}/{8}$.  Then there are constants $\La=\La(\la)>0$ and 
$K=K(\la)>0$ with the following property.
If $\rho\:Z\ra[0,\infty]$ is a Borel function,
$B=B(p,r)\sub Z$ is a ball, and $F_1,\,F_2\sub Z$ are continua
such that $F_i\cap \frac{1}{4}B\neq\emptyset$
and $F_i\setminus B\neq \emptyset$ for $i=1,2$,
then there are disjoint balls $B_i\defeq B(q_i,\la r)$
for $i=1,2$, and a path $\si\:[0,1]\ra Z$
 such that

\begin{itemize}

\item[\rm{(i)}] $q_i\in F_i$ for $i=1, 2$,
\smallskip 

\item[\rm{(ii)}] $B_i\sub \frac{7}{8}B$  and
$$\int_{B_i}\rho^Q\, d\mu \leq 80\la\int_B \rho^Q\, d\mu
$$
for $i=1,2$,
\smallskip 

\item[\rm{(iii)}] the path
 $\si$ joins $\frac{1}{4}B_1$ and $\frac{1}{4}B_2$,  has image
contained in $\La B$,  length
at most $\La r$, and $\rho$--length at most
$$
K\left( \int_{\La B}\rho^Q\, d\mu\right)^{{1}/{Q}}.
$$

\end{itemize}

\end{lemma}
\proof We can find a  subcontinuum $E_1 \sub F_1$ which 
is  contained in $ \bar B(p,\frac{3r}{8})\setminus B(p,\frac{r}{4})$  
and  joins the  sets 
$\partial B(p,\frac{3r}{8})$ and $\partial B(p,\frac{r}{4})$.
Similarly, we can find a  subcontinuum $E_2 \sub F_2$ which
is  contained in $ \bar B(p,\frac{3r}{4})\setminus B(p,\frac{5r}{8})$
and  joins the  sets
$\partial B(p,\frac{3r}{4})$ and $\partial B(p,\frac{5r}{8})$.
   Then $\diam(E_i)\geq r/8$
for $i=1,2$, and $\dist(E_1, E_2)\ge r/4$.

Applying Lemma \ref{goodpoints}
with $X=B$, the measure $\nu$ defined by
 $\nu(N)\defeq \int_N\rho^Q\, d\mu$ for a Borel set $N\sub B$,
and $Y=E_i$ for $i=1,2$,
we find $q_i\in E_i$ such that
\begin{equation}
\label{massdrop}
\int_{B(q_i,s)}\rho^Q\,d\mu\leq \frac{10s }{\diam(E_i)}
\int_B\rho^Q\,d\mu
\leq \frac{80s}{r}\int_B \rho^Q\,d\mu,
\end{equation}
for $0<s \leq {r}/{4}$.  Set $B_i\defeq B(q_i,\la r)$.
By our assumption on $Z$, we can find a path $\si$ from 
$\frac{1}{4}B_1$ to $\frac{1}{4}B_2$ with $\im(\si)\sub \La B$,
 length at most
$\La r$, and $\rho$--length at most
$$
K\left( \int_{\La B}\rho^Q\,d\mu\right)^{{1}/{Q}},
$$
where $\La =\La (\la)>0$ and $K=K(\la)>0$.

The balls $B_1$ and $B_2$ are disjoint
since $\la<{1}/{8}$ and
$$d(q_1,q_2) \ge \dist(E_1,E_2)\ge r/4.$$
Conditions (i) and (iii) are  clearly satisfied.  Condition (ii) follows from
the
facts  that $\la<{1}/{8}$ and $q_i\in \bar B(p,\frac{3r}{4})$, and from 
\eqref{massdrop}.
\qed

\medskip
\no
{\bf Proof of Proposition \ref{ballloewner}}\qua
Fix  $0<\la<{1}/{8}$ subject to the condition
$2\cdot (80\la)^{{1}/{Q}}<1$. 

Suppose  $E_1,\,E_2\sub Z$ are nondegenerate continua, and 
let $\rho\:Z\ra[0,\infty]$ be  a Borel function. 
It is enough to 
show that there is a rectifiable path $\ga$ joining $E_1$ to $E_2$
whose $\rho$--length is at most 
$$
M\left(\int_Z\rho^Q\,d\mu\right)^{{1}/{Q}}, 
$$
where $M>0$ only depends on the relative 
distance $\Delta(E_1, E_2)$ of $E_1$ and $E_2$. 
In fact, the path 
produced will have length $\lesssim \dist(E_1,E_2)$, though
this will not be used anywhere.
We will henceforth assume that $E_1$ and $E_2$ are disjoint, for otherwise
we may use a constant path mapping into $E_1\cap E_2$.

\begin{figure}[ht!]
\cl{\includegraphics[width=2.6in]{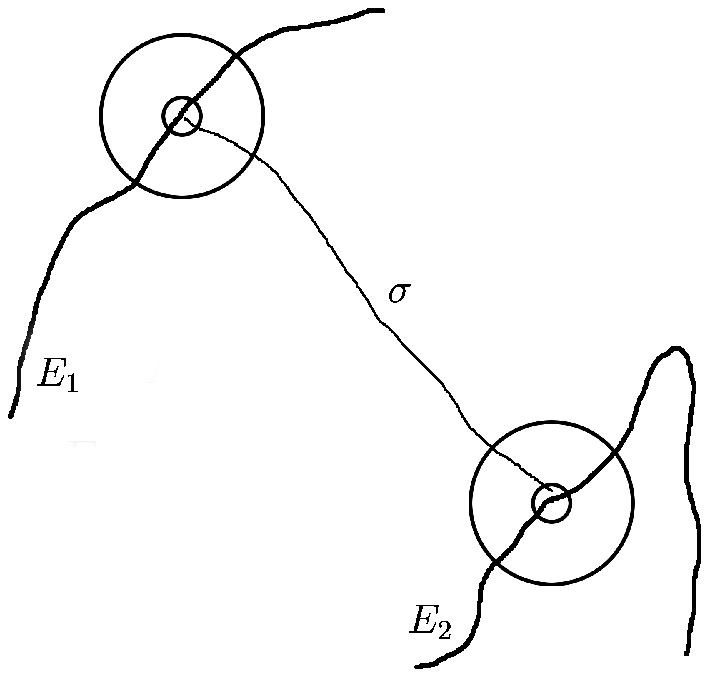}}
\nocolon\caption{}\label{fig1}
\end{figure}

Pick $p_i\in E_i$ such that $d(p_1,p_2)=\dist(E_1,E_2)$.  Set 
$$r_0\defeq \frac{1}{2}\min(d(p_1,p_2),\diam(E_1),\diam(E_2))>0.$$
Let $B_i=B(p_i,r_0)$ for  $i=1,2$.  
Then $B_1\cap B_2=\emptyset$ and $E_i\setminus B_i\neq\emptyset$ for
$i=1,2$.   Also, $\dist(B_1,B_2)\leq \left(\frac{d(p_1,p_2)}{r_0}\right)r_0
\leq tr_0$
where $t\defeq 2\max(1,\reldist(E_1,E_2))$.
By our hypotheses
we can find a path $\si$ joining $\frac{1}{4}B_1$ to
$\frac{1}{4}B_2$ of length at most $Lr_0$ and 
$\rho$--length at most
\begin{equation}
\label{sibound}
\frac{1}{m^{{1}/{Q}}}\left(\int_Z \rho^Q\,d\mu\right)^{{1}/{Q}},
\end{equation}
where $L=L(t)$, $m=m(t)$ are the constants appearing in the statement
of Proposition~\ref{ballloewner}.
Let $\Si_0\defeq\{\si\}$,  $\B_0\defeq\{ B_1, B_2\}$,
and  $\Om_0$ be the set $\{ E_1,\im(\si),E_2\}$ endowed
with the linear ordering $E_1<\im(\si)<E_2$.  In addition,
we associate the ball $B_1$ with the pair $E_1<\im(\si)$,
and the ball $B_2$ with the pair $\im(\si)<E_2$; see Figure \ref{fig1}.

Inductively, assume that for $j=0,\ldots,k$ we have
a path collection $\Si_j$, a ball collection $\B_j$,
and a collection of continua $\Om_j$ subject to the 
following conditions:

(1)\qua For $0\leq j\leq k$, the set $\Om_j$ is linearly
ordered.

(2)\qua  For each pair $\tau_1<\tau_2$ of successive elements
of $\Om_j$, there is an associated ball $B_{\tau_1,\tau_2}\in \B_j$
such that $\tau_i\setminus B_{\tau_1,\tau_2}\neq\emptyset$ 
 and $\tau_i\cap \frac{1}{4}B_{\tau_1,\tau_2}\neq\emptyset$, for
$i=1,2$.

(3)\qua For $j\geq 1$, the collections $\Si_j$, $\B_j$,
and $\Om_j$ are obtained from $\Om_{j-1}$ and $\B_{j-1}$
using the following procedure.    For each 
pair of successive elements $\tau_1,\tau_2\in\Om_{j-1}$ with
 associated ball $B_{\tau_1,\tau_2}\in\B_{j-1}$, one applies Lemma 
\ref{subballs} with $B=B_{\tau_1,\tau_2}$ and  $\{F_1,F_2\}=\{\tau_1,\tau_2\}$,
to obtain a path $\si=\si(\tau_1,\tau_2)$ and a pair
of disjoint balls 
$B_{\tau_1,\si}, B_{\si,\tau_2}$.
Here $B_{\tau_1,\si}$ is centered at a point in $\tau_1$,
and $B_{\si,\tau_2}$ is centered at a point in $\tau_2$.
Then  $\Si_j$ is the collection of these paths $\si$ and
$\B_j$ is the collection of balls $B_{\tau_1,\si},B_{\si,\tau_2}$
where $\tau_1<\tau_2$ ranges over all  successive pairs
in $\B_{j-1}$.  The continuum collection $\Om_j$
is the disjoint union $\Om_{j-1}\sqcup\{\im(\si): \si\in\Si_j\}$.
One linearly orders $\Om_j$ by extending the order on $\Om_{j-1}$
subject to $\tau_1<\im(\si)<\tau_2$; moreover, one 
associates   the ball $B_{\tau_1,\si}$
 with the  pair
$\tau_1<\im(\si)$, and the ball  $B_{\si,\tau_2}$ with
 the pair $\im(\si)<\tau_2$.

By our second induction assumption, the hypotheses of Lemma
\ref{subballs} hold for each successive pair $\tau_1,\tau_2\in\Om_k$
and associated ball $B_{\tau_1,\tau_2}\in\B_k$.
Hence we may use the procedure in the third induction
assumption (with $j$ replaced by $k+1$) to generate
$\Si_{k+1}$, $\B_{k+1}$, $\Om_{k+1}$, the linear order on $\Om_{k+1}$,
and an association of balls in $\B_{k+1}$ with successive
pairs in $\Om_{k+1}$.  The conditions in Lemma \ref{subballs}
guarantee that the induction hypotheses are fulfilled.
Therefore by induction there are collections $\Si_k$,
$\B_k$, and $\Om_k$ which satisfy
the conditions 1--3  for all $k\geq 0$. 

By induction, one proves the following:

(a)\qua For each $k\geq 0$, we have $\#\B_k=2^{k+1}$,
and each ball in $\B_k$ has radius $\la^{k}r_0$ (see Lemma \ref{subballs}).

(b)\qua  For each $k\geq 0$, $0<j\leq k$, and  $B\in \B_k$, there is
a ball $B'\in \B_{k-j}$ such that $B\sub \frac{7}{8}B'$,
(see condition (ii) of Lemma \ref{subballs}).

(c)\qua  For $k\geq 0$, the $\rho$--mass of each ball $B\in \B_k$
satisfies
$$
\int_B \rho^Q\,d\mu \leq (80\la)^k\int_Z \rho^Q\, d\mu
$$ 
(see Lemma \ref{subballs}, condition (ii)).

(d)\qua  For each $k>0$, we have $\#\Si_k=2^k$.  Each
$\si\in\Si_k$ has length at most $\La \la^{k-1} r_0$,  for
a suitable  ball $B\in\B_{k-1}$ we have $\im(\si)\sub \La B$, 
and the $\rho$--length of $\si$ is at most
$$
K\left(\int_{\La B}\rho^Q\,d\mu\right)^{{1}/{Q}},
$$
where $\Lambda$ and $K$ are as in 
Lemma~\ref{subballs}.  

(e)\qua  For $k\geq 0$ set
$$
Y_k\defeq \bigcup_{j=0}^k\bigcup_{\si\in\Si_j}\im(\si).
$$
Then $\dist(Y_k,E_1),\dist(Y_k,E_2)\leq {\la^kr_0/}{4}$,
and $Y_k$ is $({\la^kr_0}/{4})$--connected: given $y,\,y'\in Y_k$
there are points $y=y_1,\ldots,y_l=y'$ such that
$d(y_j,y_{j+1})\leq {\la^kr_0}/{4}$
for $1\leq j\leq l$.

(f)\qua  For $k>0$, the union 
$$
Z_k\defeq \bigcup_{j>k}\left(\bigcup_{\si\in\Si_j}\im(\si)\right)
$$
is covered by the collection $\{(\La +1)B\}_{B\in\B_k}$.

\bigskip
Set
$$
Y\defeq \bigcup_{j=0}^\infty\bigcup_{\si\in\Si_j}\im(\si).
$$ 
By (e) the closure $\bar Y$ of $Y$ is compact, connected, and intersects
both $E_1$ and $E_2$.   Then (f) and (a) imply that $\bar Y\setminus Y$
has $1$--dimensional Hausdorff measure zero.  Combining
this with (d) and the fact that $\la<{1}/{8}$, 
 we get that the $1$--dimensional Hausdorff
measure of $\bar Y$ is at most
$$\sum_{k=0}^\infty\sum_{\si\in\Si_k}\length(\si)\leq
Lr_0+ \sum_{k=1}^\infty(2^k)(\La\la^{k-1}r_0)
=Lr_0+\frac{2\La r_0}{1-2\la}<\infty. 
$$
Hence there is a rectifiable path $\ga\:[0,1]\ra Z$
contained in $\bar Y$ joining $E_1$ to $E_2$ with
$$\length(\ga)\leq Lr_0+ \frac{2\La r_0}{1-2\la}.$$ 
Moreover,  we may assume  that  $\im(\ga)$ is an arc and
$\ga$   is an injective map. 

Pick an  integer $s\ge 2$ such that $(1/\la)^{s-1}>8\La $. 
Then by (a) and  (b)  for every  $k\ge s$ and every 
ball $B\in \B_{k-1}$, there exists a ball $B'\in \B_{k-s}$ such that 
$\Lambda B\sub B'$. 
Hence by  (b) and (c), for every $k\geq s$ and  $\si\in \Si_k$ there is
a ball  $B'\in\B_{k-s}$ such that the $\rho$--length of $\si$
is at most
$$K\left(\int_{B'}\rho^Q\,d\mu\right)^{{1}/{Q}}
\leq K\left( (80\la)^{k-s}\int_Z\rho^Q\,d\mu\right)^{{1}/{Q}}.
$$
For $0<k<s$, each $\si\in\Si_k$ has $\rho$--length at most
$$
K\left(\int_Z\rho^Q\,d\mu\right)^{{1}/{Q}}.
$$
Recall that $\si\in\Si_0$ has $\rho$--length at most 
$$
\frac{1}{m^{{1}/{Q}}}\left(\int_Z \rho^Q\,d\mu\right)^{{1}/{Q}}.
$$
Using these bounds for $\rho$--length and the fact that $\ga$ parametrizes
an arc, we  get 
\begin{eqnarray*}
\int_\ga\rho \,ds
&\leq& \int_{\bar Y}\rho \,d\h^1 \, \leq \,  
\sum_{k=0}^\infty\sum_{\si\in\Si_k}\int_{\im(\si)}\rho \,d\h^1 \\
&\leq & \left(\int_Z \rho^Q\,d\mu\right)^{{1}/{Q}} 
\left(
\frac{1}{m^{{1}/{Q}}}
+K\sum_{k=1}^{s-1}2^k
+K\sum_{k=s}^\infty 2^k
(80\la)^{(k-s)/Q} \right) \\
&=: & M\left(\int_Z\rho^Q\, d\mu\right)^{{1}/{Q}}. 
\end{eqnarray*} 
Note that  $M$ is finite  since
$2\cdot(80\la)^{{1}/{Q}}<1$ by our initial choice of $\la$. Moreover, 
$M$ depends only  on $\Delta(E_1, E_2)$.  This shows that the 
path $\ga$ has the desired properties.
\qed

\section{Rescaling and abundance of thick paths}
\label{preplem}
\no
We now let $(Z,d,\mu)$ be an Ahlfors $Q$--regular compact  metric
space, $Q>1$,  which carries a  family of nonconstant paths with positive 
$Q$--modulus, and we let $G\acts Z$ be a uniformly
\qm\ action which is fixed point free, and acts cocompactly on triples
in $Z$. As we have seen, Lemma~\ref{aethick} implies that there exist 
nonconstant
 thick paths in $Z$.

\begin{lemma}
\label{densitypoint}
There exist disjoint open  balls $B$ and $B'$ in $Z$  such that the
set of endpoints of thick paths connecting $B$ and $B'$
has a point of density in $B$.
\end{lemma}
\proof
Let $\ga\:I\ra Z$ be a nonconstant thick path.
Pick  $t\in(0,1)$ so  that $x\defeq\ga(t)$ is distinct 
from the endpoints $x_0:=\ga(0)$ and $x_1:=\ga(1)$ of $\ga$.
Define  $R\defeq \frac{1}{10}\dist(x, \{x_0,x_1\})$,
and let $\eps\defeq R/10$. Set  $B\defeq B(x,R)$
and $B':=B(x_1,R)$. 

Every path  $\al\in B(\ga,\eps)$ has an image intersecting 
 the open
ball $B$ and picks up length in $B$ which is comparable to $R$. 
In particular, each path in $B(\ga,\eps)$ is nonconstant. 
Let 
$$
S\defeq B\cap \{\im(\al): \al\in B(\ga,\eps) \cap \P_T\}.
$$ 
By Lemma~\ref{thickproperties},  every point in 
$S$ is the initial point of a thick path ending in $B'$. 
Hence it is enough to show that
$S$ has positive $Q$--dimensional Hausdorff measure.
If this is not the case, we can find a Borel set $S'\supseteq S$ of vanishing 
Hausdorff $Q$--measure. 
Then the   function $\rho\:Z\ra[0,\infty]$ defined to be infinity
on $S'$ and $0$ elsewhere is Borel and 
  an  admissible test function
for the path  family $B(\ga,\eps)\cap \P_T$. Since the total $Q$--mass 
of $\rho$ is zero, 
Lemma \ref{aethick} implies
$$\modq(B(\ga,\eps)\setminus \C)
=\modq(B(\ga,\eps))=\modq(B(\ga,\eps)\cap \P_T)=0,  $$
which contradicts the thickness of $\ga$.
\qed

\begin{lemma}
\label{pairdensity}
Let $M\sub Z\times Z$ be the set of pairs of points
that can be joined by a thick path.  Then $M$ is dense
in $Z\times Z$.
\end{lemma}
\no
Note that this implies in particular that $Z$ is connected. 
\proof
By Lemma \ref{densitypoint}, we can find  a pair of disjoint
open balls $B$ and $B'$ so that there exists a density point  $x\in B$ 
of the set of initial points of the family  $\Ga$ 
of thick paths starting
in $B$ and ending in $B'$. For $k\in \N$ pick  $R_k>0$
with $\lim_{k\ra\infty}R_k=0$, and let
$D_k$ be the set of initial points of paths in $\Ga$
which start in $B_k\defeq B(x,R_k)$.  Then 
$$
\eps_k\defeq \frac{\hd(D_k,B_k)}{R_k}\ra 0 \quad\text{as}\quad k\to \infty.
$$
We let $\de_k\defeq\sqrt{\eps_k}$ and
define  $x^1_k\defeq x$. Since $Z$ is uniformly perfect,
for large $k$
we can choose  points $x^2_k, x^3_k\in Z$ such that the distances 
$d(x_k^1, x_k^2)$, $d(x_k^1, x_k^3)$, $d(x_k^2, x_k^3)$ exceed 
$\de_k R_k$, but are not larger than  $\la_kR_k$, 
where $\la_k=C\delta_k$ and  $C\ge 1$ is a constant independent of $k$.
Using the cocompactness
of the action $G\acts \trip(Z)$, we can  find 
$g_k\in G$ such that the image of the triple
$(x^1_k,x^2_k,x^3_k)$ under $g_k$ is $\de'$--separated
where $\de'>0$ is independent of $k$.
Applying Lemma~\ref{fillsin}, we 
conclude that 
$\hd(g_k(D_k),Z)\ra 0$
and $\diam(Z \setminus g_k(B_k))\ra 0$ as $k\ra\infty$.
After passing to a subsequence if necessary,
we may assume that the sets $Z\setminus g_k(B_k)$
Hausdorff converge to $\{z\}$ for some $z\in Z$. Since $B'$ and $B_k$ 
are disjoint for large $k$, we then also  have $\hd(g_k(B'), \{z\})\to 0$
as $k\to \infty$.

Now let 
$z_1,z_2\in Z$ and $\eps>0$ be arbitrary.   By Lemma~\ref{qmproperties},
we can find  $g\in G$ such that $g(z) \in B(z_1,\eps)$; then 
$g\circ g_k(B')\sub B(z_1, \eps)$ for large $k$. 
In addition, for large  $k$  we will
also have $g\circ g_k(D_k)\cap B(z_2,\eps)\neq\emptyset$. 
Using this and the invariance of thickness under \qm\ homeomorphisms
we see that  there is a thick path starting in 
$B(z_1,\eps)$ and ending in $B(z_2,\eps)$.
\qed

\begin{lemma}
\label{conctballs}
For each  $C>0$ there are constants $m>0$ and 
$L>0$ such that if $B,\,B'\sub Z$ 
are $R$--balls with $\dist(B,B')\leq CR$,
then the modulus of the family of paths
of length at most $LR$ joining $B$ to $B'$
has $Q$--modulus greater than  $m$.

\end{lemma}
\proof Suppose $C>0$ is arbitrary. 
We first claim that  there is a number  $m_0>0$
such that if $B,\,B'\sub Z$ 
are $R$--balls with $\dist(B,B')\leq CR$, then
$\modq(\Ga(B,B'))>m_0$.
If this assertion were false, there would be  
balls  $B_k=B(z_k,R_k)$ and $B_k'=B(z_k',R_k)$ for $k\in \N$
such that $\dist(B_k,B_k')\leq CR_k$ for all $k$, 
but 
\begin{equation}
\label{tendstozero}
\lim_{k\ra\infty}\modq(\Ga(B_k,B_k'))=0.
\end{equation}
Passing to a subsequence, we may assume that the
sequences $(z_k)$ and $(z_k')$ converge  in $Z$. 
Lemma \ref{pairdensity} then implies 
that $\lim_{k\ra\infty}R_k=0$.   Disregarding finitely many sequence elements
if necessary, 
 we can choose
$x^1_k\in \partial B(z_k,R_k)$, and $x^2_k\in \partial B(z_k,2R_k)$ for all $k$
by the connectedness of $Z$. 
Pick  $g_k\in G$ such that the triples $(g_k(z_k),g_k(x^1_k),g_k(x^2_k))$
are $\de$--separated where $\de>0$ is independent of $k$.
Since the homeomorphisms
$g_k$ are uniformly \qm\ and the balls $B_k$ and $B_k'$ 
have uniformly bounded relative distance, it is easy to see that there is 
$\eps>0$ such that $B(g_k(z_k),2\eps)\sub g_k(B_k)$
and $B(g_k(z_k'),2\eps)\sub g_k(B_k')$ for all $k$.
Passing to yet another subsequence, we may assume that the sequences 
$(g_k(z_k))$ and $(g_k(z_k'))$ converge to points $z\in Z$ and $z'\in Z$,
respectively. 
Hence for large $k$ we have $\tilde B:=B(z,\eps)\sub g_k(B_k)$,
$\tilde B':=B(z',\eps)\sub g_k(B_k')$. Tyson's theorem
(Theorem \ref{qmodpreserved}) gives
$$\modq(\Ga(B_k,B_k'))\geq c \modq\big(\Ga(g_k(B_k),g_k(B_k'))\big)
\geq c\modq(\Ga(\tilde B, \tilde B')),
$$
where $c>0$ is a constant independent of $k$.
Since $\modq(\Ga(\tilde B, \tilde B'))>0$ by Lemma
\ref{pairdensity}, this contradicts (\ref{tendstozero}), and hence
the claim is true. 

According to Lemma~\ref{qmodproperties} we can choose $L\ge 2$ 
such that every  family of paths in $Z$  which start 
in a given 
ball of radius $R$ and have length at least $LR$ has modulus at most 
$m_0/2$. Now  $B$ and $B'$ are  arbitrary 
balls 
of radius $R>0$ in $Z$ 
 and let $\Gamma_1$ and $\Gamma_2$ be the families  of paths in 
$Z$ a which connect $B$ and $B'$ and have length at most $LR$ 
and length at least $LR$, respectively. 
Then by the choice of $L$ and by  subadditivity of modulus we have 
$$ m_0 < \modq(\Ga(B,B'))\le \modq(\Gamma_1)+\modq(\Gamma_2) 
\le \modq(\Gamma_1) +m_0/2.$$
It follows that $\modq(\Gamma_1) > m:= m_0/2>0$.   
\qed

\section{The proofs of the theorems}
\label{proofsofthms}

 \no
{\bf Proof of Corollary \ref{specialkeithlaakso}}\qua
Under the  hypotheses of the corollary,  for every weak
tangent $W$ of $Z$ there exist a point $z\in Z$ and a \qm\ homeomorphism
between $W$ and $Z\setminus \{z\}$ (cf\ \cite[Lemma~5.2]{quasimobius}). 
According to Theorem~\ref{keilaa} there exists  a weak tangent $W$ of 
$Z$ which carries a family of nonconstant paths with positive $Q$--modulus. 
Note that $W$ is an Ahlfors $Q$--regular complete metric space. 
Hence by Tyson's Theorem~\ref{qmodpreserved}, the space $Z$ also carries a 
family of nonconstant paths with positive $Q$--modulus. 
\qed 

\medskip \no
{\bf Proof of Theorem \ref{loewner}}\qua   If $Z$ is as in the theorem, 
then $Z$ satisfies the assumptions made in the begining 
of Section~4.  So  Lemma \ref{conctballs} applies, showing  that 
$Z$ satisfies the hypotheses of
Proposition~\ref{ballloewner}.
Therefore, $Z$ is a $Q$--Loewner space.
\qed

\medskip \no
{\bf Proof of Theorem \ref{mainthm}}\qua  Let $G$ be as in the statement
of Theorem \ref{mainthm},  suppose $Z$ is an Ahlfors $Q$--regular metric
space where $Q\ge 2$ is the Ahlfors regular
conformal dimension of $\geo G$, and $\phi\:\geo G\ra Z$
is a \qs\ homeomorphism.  Conjugating the canonical action 
$G\acts\geo G$ by $\phi$, we obtain a uniformly \qm\ action
$G\acts Z$ which is fixed point free and for which the
induced action on triples is both properly discontinuous and cocompact.
Now by Corollary~\ref{specialkeithlaakso} the space $Z$ 
carries a 
family of nonconstant paths with positive $Q$--modulus, and so it is
 $Q$--Loewner by Theorem \ref{loewner}. 

 According to
 Theorem 1.2 of \cite{qparametr} every Ahlfors  $Q$--regular
and $Q$--Loewner $2$--sphere is quasisymmetrically  homeomorphic 
to the standard $2$--sphere $\sph$. This applies to $Z$ and so there 
exists a  \qs\ homeomorphism $\psi\:Z\ra \sph$. 
Conjugating our action
$G\acts Z$ by $\psi$, we get a uniformly quasiconformal
action $G\stackrel{1}\acts\sph$ (we use the superscript ``1"
to distinguish this action from another  action discussed below). 
 By a theorem of Sullivan and Tukia
(cf\ \cite[p.~468]{sul}
and \cite[Theorem~F and Remark~F2]{Tuk})  every   uniformly
quasiconformal action on $\sph$  is quasiconformally conjugate to an action 
by   M\"obius transformations. Hence $G\stackrel{1}\acts\sph$
is quasiconformally conjugate to an action
$G\stackrel{2}{\acts} \sph$ by M\"obius transformations.
If we represent $\H^3$ by the unit ball model so  
that $\geo \H^3=\sph$, the action $G\stackrel{2}{\acts} \sph$
naturally extends to an isometric  action $G\acts\H^3$. 
Being topologically conjugate to  $G\stackrel{1}\acts\sph$ and hence to 
$G\acts Z$, the action 
$G\stackrel{2}\acts\sph$ is also  properly discontinuous
and cocompact on triples.
Therefore, the corresponding
isometric action  $G\acts \H^3$ is discrete and cocompact.
\qed 

\medskip \no {\bf Proof of Theorem~\ref{kleinian}}\qua 
Let  $G$ be a group  as  in the theorem, and assume that $Q>1$ 
is equal to the Ahlfors regular conformal dimension of $Z=\Lambda(G)$. 
Note that $Z\sub \mathbb{S}^n$ 
equipped with the ambient Euclidean  metric on $\mathbb{S}^n\sub 
\mathbb{R}^{n+1}$ is Ahlfors
$Q$--regular \cite[Theorem~7]{sul2},
and that the induced action $G\acts Z$ satisfies the hypotheses 
of Theorem~\ref{loewner}. It follows that $Z$ is a $Q$--Loewner 
space. Hence $Z$ satisfies a $(1,Q)$--Poincar\'e inequality. 

Since the metric space   $Z\sub \R^{n+1}$ is Ahlfors
regular and satisfies a $(1,Q)$--Poincar\'e inequality,    
a theorem by Cheeger \cite[Theorem\ 14.3]{cheeger} implies that $Z$ has 
a weak tangent which is bi-Lipschitz equivalent to some Euclidean 
space  $\R^k$, $k\ge 1$. As in the proof of Corollary~1.6, every weak 
tangent of $Z$ is homeomorphic to $Z$ minus a point. 
Hence $Z$ has topological dimension $k$. 

On the other hand, 
since   each  weak tangent of 
an Ahlfors  $Q$--regular space is also Ahlfors  $Q$--regular, we conclude 
that $Q=k\in \N$;  in particular,   the topological 
dimension of $Z$ is equal to its Hausdorff dimension. 
Moreover, since $Q>1$, we have $k\ge 2$.  
The desired conclusion now follows from \cite[Theorem\ 1.2]{quasimobius}.
\qed 

\medskip \no
The proof shows that if $G\acts X$
is a properly discontinuous, quasi-convex cocompact, and
isometric action on a $\mathrm{CAT}(-1)$-space $X$,
if the limit set $\Lambda(G)$ has Hausdorff dimension equal to its Ahlfors
regular conformal dimension, and if   $\Lambda(G)$ 
embeds in some Euclidean space by a bi-Lipschitz map, 
then there is a convex $G$--invariant
copy of a  hyperbolic space $Y\sub X$
 on which $G$ acts cocompactly.

The conclusion of   the  quoted result by  Cheeger already holds
if $Z$ satisfies a $(1,p)$--Poincar\'e inequality for some $p>1$. 
So it would be enough to stipulate this condition in Theorem~\ref{kleinian}
instead of requiring that the Ahlfors regular conformal dimension
of $Z=\Lambda(G)$ is equal to $Q$. 

The converse of Theorem~\ref{kleinian}  and its indicated modifications
lead to a statement that is  worth
recording: 
 if the limit set $\Lambda(G)$ of a group
 $G$ as in the theorem  is {\em not}   a ``round"
subsphere of $\mathbb{S}^n$,
then the Ahlfors regular dimension of $Z=\Lambda(G)$ is strictly
less than $Q$, $Z$ does not  carry a family of
nonconstant curves  with positive $Q$--modulus, 
and
  $Z$ does not  satisfy
a $(1,p)$--Poincar\'e inequality for any $p>1$.
In particular, limits sets of such groups  $G$ cannot lead to new examples 
of Loewner spaces.

\section{Spaces whose Ahlfors regular conformal dimension is not realized}
\label{nonrealized}

\no In our discussion below,
we will refer to the Ahlfors regular conformal dimension
simply as the conformal dimension. 

The most basic  example of a space whose  conformal
dimension is not realized is the standard Cantor set $C$.  
This dimension is equal to $0$ for $C$, but it is not attained, since 
any  quasisymmetric homeomorphism between $C$ and a metric space 
$Z$  is  bi-H\"older, 
and so the   Hausdorff dimension of $Z$ is strictly positive. 

To our knowledge, the 
first connected and locally 
 connected example of this type is due to Pansu, which
we learned of through M.~Bourdon.  Essentially the same example
was considered also in \cite{buyalo}:
if one glues together two closed hyperbolic surfaces $N_1$ and
$N_2$  isometrically along
embedded  geodesics $\ga_i\sub N_i$ of equal length, then 
one  obtains  a $2$--complex $M$ with
curvature bounded from above by $-1$ and 
the boundary at infinity $\geo\widetilde M$
of the universal cover $\widetilde M$
has  conformal dimension $1$. 
To see  this one  pinches  the hyperbolic structures along the
closed geodesics $\ga_i$, and observes that the volume entropies
of the resulting universal covers tend to $1$ (``branching becomes less and
less frequent''). The space $\geo \widetilde M$ is not quasisymmetrically
homeomorphic to
an Ahlfors $1$--regular space, because in this case it would have 
to be a topological circle by \cite[Theorem 1.1]{quasimobius};
in fact it is not difficult to
see directly that
$\geo\widetilde M$ is not quasisymmetrically  homeomorphic to 
a space with finite
$1$--dimensional Hausdorff measure.

Bishop and Tyson \cite{bistys}
have shown that ``antenna sets''---certain
self-similar dendrites in the plane---have conformal dimension $1$,
but are not quasisymmetrically homeomorphic to any space of 
Hausdorff dimension $1$.

Another example of a similar flavor is due to Laakso. He has    
shown that the
standard Sierpinski gasket has  conformal
dimension $1$, but again, this dimension  cannot be realized.
By considering pairs of points whose removal disconnects the  set, 
it is not hard to show  that
 the homeomorphism group of the gasket is the same as its
isometry group for the usual embedding in $\R^2$. It follows that 
 this example
is not homeomorphic to Pansu's example.

There are  translation invariant  Ahlfors
regular metrics on $\R^2$
for which the $1$--parameter group of linear
transformations $e^{tA}$,
where
$$
A\defeq \left[\begin{matrix}
1&1\\
0&1\\
\end{matrix}
\right],
$$
is a family of homotheties.  Their  conformal
dimension is $2$, but it cannot be  realized. 
The second author would like to thank
L.~Mosher for bringing these examples to his attention.
One can also describe them  as follows.
Let $G$ be the semi-direct product of 
$\R$ with $\R^2$, where $\R$ acts on $\R^2$ by the $1$--parameter
group above.  Then the solvable Lie group $G$ admits left
invariant Riemannian metrics with curvature pinched arbitrarily close
to $-1$; if one removes the unique fixed point from $\geo G$,
one gets the ``twisted plane'' example above.  

The examples discussed so far  are all either disconnected,  have
local cut points (ie, by removing a single point one
can disconnect a connected neighborhood), or cannot be the boundary
of a hyperbolic group.

  Suppose an Ahlfors  $Q$--regular space $Z$
is quasisymmetric to the boundary of a hyperbolic group $G$,
where $Q$ is the  conformal dimension
of $Z$.  If $Q<1$, then the topological dimension of $Z$ is zero;
thus there is a free subgroup $F_k$ sitting in $G$ with 
finite  index.   But then $k=1$,  $|Z|=2$, and $Q=0$, for
otherwise $k>1$ which  implies  that
$Z$ is quasisymmetric to the standard Cantor set, whose
conformal dimension is not realized.  If $Q=1$, then 
\cite[Theorem 1.1]{quasimobius} implies that $Z$ is 
quasisymmetric to the standard circle.
The case $Q>1$ is covered by Theorem \ref{loewner}. Disconnected spaces,
or spaces with local
cut points cannot satisfy a  
Poincar\'e inequality, so Theorem \ref{loewner} implies
that  if  $Q>1$, then  $Z$ is connected and has
no local cut points.   The examples of Bourdon and Pajot
\cite{bp3} give
boundaries of hyperbolic groups which
possess  these two topological properties,
but which are not quasisymmetrically
homeomorphic to a $Q$--regular space satisfying a 
$(1,Q)$--Poincar\'e inequality.  Thus by Theorem 
\ref{loewner} even under these  topological
conditions the conformal dimension is not 
 necessarily realized. 

Based on the examples mentioned above, one may speculate that if the
 conformal dimension of a self-similar
space fails to be attained, then this
is due to degeneration which leads to a limiting structure
resembling a foliation or lamination.

We conclude this section with two questions related to the
realization problem. 

\begin{problem}
Can one  algebraically characterize 
the hyperbolic groups whose
boundary has (Ahlfors regular) conformal dimension equal to $1$? In particular, 
if the boundary of a hyperbolic group is homeomorphic to
a Sierpinski carpet or a Menger curve, is
the Ahlfors regular conformal dimension strictly greater than $1$?
\end{problem}

\begin{problem}
Is the (Ahlfors regular) conformal dimension of the standard
square Sierpinski carpet $S$ attained?   
\end{problem}

\noindent
If it is, 
it  seems  to be the case  that $S$ equipped with  any Ahlfors regular metric
realizing its   conformal dimension is  a Loewner
space.
We remark that it follows from \cite{keilaa} that the 
 conformal dimension of $S$ is less than
its Hausdorff dimension.  A calculation by the second
author had earlier given an explicit  upper bound for the conformal
dimension of $S$.

\section{Remarks and open problems}
\label{problems}

\no 
The themes explored in this paper lead to various general questions.  
To further exploit the relation between the algebraic structure of 
a Gromov hyperbolic group and the analysis of its boundary one needs
 analytic tools from the general theory of analysis
on metric spaces, perhaps tailored to the setting of self-similar spaces or 
spaces admitting group actions as considered in this paper. 
In particular, it would be interesting to find  classes of function
spaces that are invariant under quasisymmetric homeomorphisms. 
They could be used to define quasisymmetric invariants 
and answer structure and rigidity questions for quasisymmetric homeomorphisms.

The setting of Loewner spaces is relatively well-understood, but it
is not clear how natural this framework really  is. At present there is a 
somewhat limited supply of these spaces, and one would like to have 
more examples. As Theorem~\ref{kleinian} indicates, 
the Loewner condition   seems to lead to strong 
conclusions in the presence of group actions and probably also
in the presence of  self-similarity.
In view of this theorem the following problem suggests itself.   

\begin{problem}
Can one classify all quasi-convex cocompact isometric actions 
$G \acts X$, where  
$X$ is a  $\mathrm{CAT}(-1)$--space  
and the Ahlfors regular conformal dimension   
of the limit set $\Lambda(G)$ is realized and strictly 
greater than $1$?
\end{problem}

\no 
Note that in this situation $Z=\Lambda(G)$ is a Loewner space, so the problem
asks for a classification of all Loewner spaces that arise as  limit sets
of quasi-convex cocompact isometric group actions on
$\mathrm{CAT}(-1)$--spaces. 

Conversely, one could start with an Ahlfors $Q$--regular 
 $Q$--Loewner $Z$ space quasisymmetrically 
homeomorphic to the boundary of  a Gromov hyperbolic group $G$ 
and ask whether $Z$ appears  as the limit set  $\Lambda(G)$
of some isometric  action $G\acts X$, where $X$ is a negatively 
curved metric space. It is natural to require
that $X$ is Gromov hyperbolic. One can interpret the relation  between 
$Z$ and $\Lambda(G)$ in a measure theoretic sense. 
The obvious  measure on $Z$ is Hausdorff $Q$--measure, and the 
 measure on $\Lambda(G)$ related to the dynamics $G\acts 
\Lambda(G)$ is the so-called  Patterson--Sullivan measure (cf~\cite{coo}).
We arrive at the following question: 

\begin{problem}
Suppose $\phi\:Z\ra \geo G$ is a quasisymmetric homeomorphism
from a compact Ahlfors regular Loewner space $Z$
to the boundary $\geo G$ of Gromov hyperbolic group $G$.
  Is there a discrete, cocompact,
isometric action
$G\acts X$ of $G$ on a Gromov hyperbolic space
$X$ whose Patterson--Sullivan measure 
lies the same measure class as push-forward of
Hausdorff measure under $\phi$?
\end{problem}

\no
More generally, one may ask when the measure class 
of a given  measure  on the
boundary of a Gromov hyperbolic group is represented
by the  Patterson--Sullivan measure for some
Gromov hyperbolic ``filling" $G\acts X$ of the boundary action
$G\acts\geo G$.

The general problem   behind Cannon's conjecture is the desire 
to find canonical metric spaces on which a given Gromov  hyperbolic 
group $G$ acts. Since the dynamics of an isometric action 
$G$ on a Gromov hyperbolic space $X$  is encoded in the
 Patterson--Sullivan measure on $\Lambda(G)$, a first step in this direction is
to find a natural measure, or at least a natural measure class on 
$\geo G$.  

\begin{problem} \label{measureclass}
Given a Gromov hyperbolic group $G$, 
when is there a natural 
measure class   on $\geo G$? 
\end{problem}

\no
Here ``natural'' can be interpreted in various ways.
One could require the measure class  to be invariant under
all (local) quasisymmetric homeomorphisms.  For instance, if
$G$ acts discretely cocompactly on a rank $1$
symmetric space $X$  other than $\H^2$,
then the measure class of the Lebesgue  measure associated with the
standard smooth structure on $\geo X$ is invariant
under the full group of quasisymmetric self-homeomorphisms
of $\geo X\simeq \geo G$.  When $X=\H^2$ this fails,
since  the ``Mostow map'' $\geo \H^2\ra\geo \H^2$ induced between
two non-conjugate discrete,  cocompact and  isometric actions
$G\stackrel{1}{\acts} \H^2$ and $G\stackrel{2}{\acts}\H^2$
will not be absolutely continuous with respect to
Lebesgue measure.  One expects a similar phenomenon whenever
$G$ virtually splits over a virtually cyclic group, or 
equivalently, when 
$\geo G$ has local cut points.
Due to this,  one can hope for an
affirmative answer to Problem~\ref{measureclass} only
under the assumption that $G$ does not have this property.

In many cases one expects that $G$ is a subgroup of finite index 
in the group $\mathrm{QS}(\geo G)$
of all quasisymmetric self-homeomorphisms  of $\geo G$.
Then  the requirement that the measure class be 
invariant under $\mathrm{QS}(\geo G)$ is rather weak. 
This suggests another
 (stronger)  interpretation of Problem~\ref{measureclass}:
 the measure class  should be constructed in a quasisymmetrically
invariant way.


\begin{thebibliography}

\bibitem{bistys}
\textbf{Christopher~J Bishop}, \textbf{Jeremy~T Tyson}, \emph{Conformal
  dimension of the antenna set}, Proc. Amer. Math. Soc. 129 (2001) 3631--3636
  \MR{1860497}

\bibitem{qparametr}
\textbf{Mario Bonk}, \textbf{Bruce Kleiner}, \emph{Quasisymmetric
  parametrizations of two-dimensional metric spheres}, Invent. Math. 150 (2002)
  127--183
  \MR{1930885}

\bibitem{quasimobius}
\textbf{Mario Bonk}, \textbf{Bruce Kleiner}, \emph{Rigidity for
  quasi-{M}\"obius group actions}, J. Differential Geom. 61 (2002) 81--106
  \MR{1949785}

\bibitem{BKR}
\textbf{Mario Bonk}, \textbf{Pekka Koskela}, \textbf{Steffen Rohde},
  \emph{Conformal metrics on the unit ball in {E}uclidean space}, Proc. London
  Math. Soc. 77 (1998) 635--664
  \MR{1643421}

\bibitem{bourdonpajot}
\textbf{Marc Bourdon}, \textbf{Herv{\'e} Pajot}, \emph{Rigidity of
  quasi-isometries for some hyperbolic buildings}, Comment. Math. Helv. 75
  (2000) 701--736
  \MR{1789183}

\bibitem{bp3}
\textbf{Marc Bourdon}, \textbf{Herv{\'e} Pajot}, \emph{Cohomologie $\ell_p$ et
  espaces de {B}esov}, J. Reine Angew. Math. 558 (2003) 85--108
  \MR{1979183}

\bibitem{buyalo}
\textbf{Sergei Buyalo}, \emph{Volume entropy of hyperbolic graph surfaces},
  preprint available at:
  \url{http://www.pdmi.ras.ru/preprint/2000/}

\bibitem{cheeger}
\textbf{J Cheeger}, \emph{Differentiability of {L}ipschitz functions on metric
  measure spaces}, Geom. Funct. Anal. 9 (1999) 428--517
  \MR{1708448}

\bibitem{coo}
\textbf{Michel Coornaert}, \emph{Mesures de {P}atterson--{S}ullivan sur le bord
  d'un espace hyperbolique au sens de {G}romov}, Pacific J. Math. 159 (1993)
  241--270
  \MR{1214072}

\bibitem{hajlaszkoskela}
\textbf{Piotr Haj{\l}asz}, \textbf{Pekka Koskela}, \emph{Sobolev met
  {P}oincar\'e}, Mem. Amer. Math. Soc. 145 (2000) no. 668
  \MR{1683160}

\bibitem{hei}
\textbf{Juha Heinonen}, \emph{Lectures on analysis on metric spaces},
  Universitext, Springer--Verlag, New York (2001)
  \MR{1800917}

\bibitem{heikos}
\textbf{Juha Heinonen}, \textbf{Pekka Koskela}, \emph{Quasiconformal maps in
  metric spaces with controlled geometry}, Acta Math. 181 (1998) 1--61
  \MR{1654771}

\bibitem{keilaa}
\textbf{Stephen Keith}, \textbf{Tomi Laakso}, \emph{Conformal Assouad dimension
  and modulus}, {G}eom. Funct. Analysis, to appear.

\bibitem{kinshan}
\textbf{Juha Kinnunen}, \textbf{Nageswari Shanmugalingam}, \emph{Regularity of
  quasi-minimizers on metric spaces}, Manuscripta Math. 105 (2001) 401--423
  \MR{1856619}

\bibitem{laakso}
\textbf{Tomi Laakso}, \emph{Ahlfors {$Q$}--regular spaces with arbitrary {$Q>1$}
  admitting weak {P}oincar\'e inequality}, Geom. Funct. Anal. 10 (2000)
  111--123
  \MR{1748917}

\bibitem{Pau}
\textbf{Fr{\'e}d{\'e}ric Paulin}, \emph{Un groupe hyperbolique est
  d\'etermin\'e par son bord}, J. London Math. Soc. 54 (1996) 50--74
  \MR{1395067}

\bibitem{semmesquant}
\textbf{Stephen Semmes}, \emph{Finding curves on general spaces through
  quantitative topology, with applications to {S}obolev and {P}oincar\'e
  inequalities}, Selecta Math. (N.S.) 2 (1996) 155--295
  \MR{1414889}

\bibitem{semmesnovel}
\textbf{Stephen Semmes}, \emph{Some novel types of fractal geometry}, Oxford
  Mathematical Monographs, The Clarendon Press, Oxford University Press, New
  York (2001)
  \MR{1815356}

\bibitem{sul2}
\textbf{Dennis Sullivan}, \emph{The density at infinity of a discrete group of
  hyperbolic motions}, Inst. Hautes \'Etudes Sci. Publ. Math. 50 (1979)
  171--202
  \MR{0556586}

\bibitem{sul}
\textbf{Dennis Sullivan}, \emph{On the ergodic theory at infinity of an
  arbitrary discrete group of hyperbolic motions}, from: ``Riemann surfaces and
  related topics (Stony Brook, 1978)'', Princeton Univ. Press, Princeton, N.J.
  (1981)  465--496
  \MR{0624833}

\bibitem{Tuk}
\textbf{Pekka Tukia}, \emph{On quasiconformal groups}, J. Analyse Math. 46
  (1986) 318--346
  \MR{0861709}

\bibitem{tyson}
\textbf{Jeremy~T Tyson}, \emph{Quasiconformality and quasisymmetry in metric
  measure spaces}, Ann. Acad. Sci. Fenn. Math. 23 (1998) 525--548
  \MR{1642158}

\bibitem{tyson2}
\textbf{Jeremy~T Tyson}, \emph{Metric and geometric quasiconformality in
  {A}hlfors regular {L}oewner spaces}, Conform. Geom. Dyn. 5 (2001) 21--73
  \MR{1872156}

\end{thebibliography}
\end{document}